%% file: main.tex
\documentclass[
onefignum,onetabnum]{siamart171218}
\input{splitted/preamble}

\input{splitted/shared}

\ifpdf
\hypersetup{
  pdftitle={Smoothers with localized residual computations for geometric multigrid methods},
  pdfauthor={M. Wichrowski, P. Munch, M. Kronbichler, and G. Kanschat}
}
\fi

\begin{document}

\maketitle

\begin{abstract}
\input{splitted/abstract}
\end{abstract}



\input{splitted/main_text}

\bibliographystyle{siamplain}
\bibliography{literature}
\end{document}

%% file: splitted/preamble.tex
\usepackage{amstext}
\usepackage{stmaryrd}
\usepackage{algpseudocode}

\usepackage{cancel}
\usepackage{graphicx}

\usepackage{array}
\usepackage{verbatim}
\usepackage{booktabs}
\usepackage{calc}
\usepackage{textcomp}
\usepackage{multirow}

\usepackage[all]{xy}

\usepackage{mathtools}
\usepackage{placeins}

\usepackage{tikz,tikzscale}
\usetikzlibrary{calc}
\usetikzlibrary{shadings}

\usepackage[mode=multiuser,status=draft]{fixme} 
\fxsetup{innerlayout=inline}
\fxsetup{targetlayout=colorcb}

\fxusetheme{color}
\FXRegisterAuthor{cc}{envcc}{CC}
\FXRegisterAuthor{gk}{envgk}{GK}
\FXRegisterAuthor{mk}{envmk}{MK}
\FXRegisterAuthor{pm}{envpm}{PM}
\FXRegisterAuthor{mw}{envmw}{MW}

\definecolor{apply}{rgb}{0.3,.7,0.}
\definecolor{applyface}{rgb}{.9,.95,0.}
\definecolor{invert}{rgb} {0.5,0.,1.}
\colorlet{invertface}{invert!60!red}

\def\mesh{\mathbb M}

\def\R{\mathbb R}
\def\ncol{N_{\text{colors}}}
\def\rest#1{I_{#1}^\downarrow}
\def\prol#1{I_{#1}^\uparrow}

\newcommand{\pluseq}{\stackrel{+}{\gets}}

\newcommand{\LUpNew}{local\_update}
\newcommand{\LSolNew}{local\_solve}

%% file: splitted/shared.tex
\usepackage{lipsum}
\usepackage{amsfonts}
\usepackage{graphicx}
\usepackage{epstopdf}
\usepackage{amsmath,amssymb}
\ifpdf
  \DeclareGraphicsExtensions{.eps,.pdf,.png,.jpg}
\else
  \DeclareGraphicsExtensions{.eps}
\fi

\newsiamremark{remark}{Remark}
\newsiamremark{hypothesis}{Hypothesis}
\crefname{hypothesis}{Hypothesis}{Hypotheses}
\newsiamthm{claim}{Claim}

\headers{Smoothers with localized residuals}{M. Wichrowski, P. Munch, M. Kronbichler and G. Kanschat}

\title{Smoothers with localized residual computations for geometric multigrid methods\thanks{Submitted to the editors 15.12.2023.
\funding{PM and MK were partially supported by the Bayerische Kompetenznetzwerk für
Technisch-Wissenschaftliches Hoch-und Höchstleistungsrechnen
(KONWIHR) through the project ``High-order matrix-free ﬁnite
element implementations with hybrid parallelization and improved
data locality.'' GK was supported by the Deutsche Forschungsgemeinschaft (DFG, German Research Foundation) under Germany's Excellence Strategy EXC 2181/1 - 390900948 (the Heidelberg STRUCTURES Excellence Cluster)}}}

\author{MICHAŁ Wichrowski\thanks{Interdisciplinary Center for Scientifc Computing, Heidelberg University, Heidelberg, Germany}
\and PETER Munch\thanks{High-Performance Scientific Computing, University of Augsburg, Germany}
\and MARTIN Kronbichler\footnotemark[3] \thanks{Applied Numerics, Faculty of Mathematics, Ruhr University Bochum, Germany}
\and GUIDO Kanschat\footnotemark[2]
}

\usepackage{amsopn}

\makeatletter
\newcommand*{\addFileDependency}[1]{
  \typeout{(#1)}
  \@addtofilelist{#1}
  \IfFileExists{#1}{}{\typeout{No file #1.}}
}
\makeatother


%% file: splitted/abstract.tex
    We improve the performance of multigrid solvers
        on many-core architectures with cache hierarchies
    by reorganizing operations in the smoothing step to minimize memory transfers. We focus on patch smoothers, which offer robust convergence rates with respect to the finite element degree for various equations,
    in the setting of multiplicative subspace correction for numerical efficiency.
    By combining  the computation of local residuals with local solvers, we increase the locality of the problem and thus reduce data transfers.
    The thread-parallel implementation of this algorithm is based on coloring, which contradicts cache efficiency.
    We improve data locality by rearranging the loop into batches so that more data can be reused.
    The organization of consecutive batches prioritizes data locality. 

%% file: splitted/main_text.tex
\section{Introduction}

When improving the efficiency of computational algorithms, one encounters two limitations: core-bound, where performance is limited by the speed at which arithmetic computations and related operations within the compute core can be performed, and memory-bound, where performance is limited by the rate at which data can be accessed from global memory. Memory-bound algorithms can be accelerated by using caches in an efficient way.
Here, we implement a multiplicative vertex-patch smoother for multigrid methods in such a way that access to main memory is reduced. To this end, we propose an implementation which is matrix-free as in~\cite{Kormann2011}, such that only vector data must be loaded. Second, we rearrange operations such that data in caches can be reused without causing races between different threads.

It is a particular challenge of finite element codes that they are usually formulated to perform sweeps over the entire vector over and over again, which in combination with hardware capabilities results in many computational kernels in PDE solvers being limited by memory bandwidth. This is also the case for the multigrid method, which combines simple localized iteration schemes, called smoothers, on a hierarchy of meshes to
damp all error modes effectively.
Multigrid methods are among the most efficient iterative schemes for solving the systems of linear equations associated with the discretization of elliptic partial differential equations, such as Poisson's equation. A notable milestone is textbook multigrid efficiency, i.e., the solution of a linear system with a computational effort of few matrix-vector products~\cite{thomas2003textbook, kohl2022textbook}.  
An example of near-textbook efficiency is highlighted in \cite{kohl2022textbook}, where the authors successfully solve the Stokes equations on hierarchical hybrid grids, but using primarily low-order finite element methods. Achieving textbook efficiency becomes increasingly challenging for higher-order discretizations, as the number of iterations required grows with polynomial order for simple smoothers, see e.g.~\cite{Kanschat08smoother,fehn2020hybrid}.
Patch smoothers, introduced in~\cite{ArnoldFalkWinther97Hdiv,ArnoldFalkWinther00}, can effectively keep the number of iterations low regardless of the polynomial order of finite elements~\cite{WitteArndtKanschat21}. They have been proven effective in various applications such as the Stokes equations~\cite{hong2016robust,KanschatMao15}.  

Efficient application of the overlapping patch smoother requires a multiplicative approach; otherwise, the convergence rate may significantly decrease, impacting the overall effectiveness of the solver~\cite{WitteArndtKanschat21}. 
In this paper, we first notice that the sequential application of a patch smoother, as presented in~\cite{WitteArndtKanschat21},  breaks down into two parts: computation of the residual, which is a global operation through a matrix-vector product, and the application of an approximation of the inverse on a patch of small size.  We then reorganize the smoothing step in a Gauss--Seidel manner by performing matrix-vector products on parts of the vector in close temporal proximity to the inverse. 
This combination exploits the inherent locality of the problem and reduces transfers from main memory, since vector entries loaded for the matrix-vector product can be hit in cache when applying the patch inverse.  We seek further improvements by arranging the loop over patches so that the data fetched for one patch could be partially reused for the next one. This is obtained by processing patches following a space-filling curve.

The idea of performing multiple algorithmic steps during a single sweep over all spatial points has been used by the finite difference community through
temporal wave-front blocking
techniques, see, e.g., \cite{malas2017multidimensional}. 
For higher-order
finite element methods, the dependency region of the adjacent elements increases when measured in terms of the number of degrees of freedom, especially in the common element-by-element evaluation approach, which diminishes the benefit of caches.
Therefore, previous work on data locality in the finite element context has considered either
the data locality within the operator evaluation considered in isolation~\cite{Kronbichler2017a}, or the vector operations surrounding a single operator evaluation, such as embedding vector operations into the loop of matrix-free evaluation of the matrix action for point Jacobi/Chebyshev smoothers \cite{kronbichler2019multigrid,munch2023cache}, combining vector updates and inner products with the matrix-vector product
\cite{kronbichler2022cg}, flux terms in discontinuous Galerkin methods \cite{Trojak2022}, or predictor-corrector time discretizations \cite{Dumbser2018} to the best of the authors' knowledge.

Although the above concept works well in sequential implementations, it has no straightforward extension to the parallel setting. In this paper, we focus on shared memory computations, where data races in the parallel algorithm need to be avoided by tracking data dependencies. 
A commonly employed solution is coloring, a technique that categorizes patches into distinct, non-overlapping sets for parallel processing. While this minimizes the need for synchronization, allowing for efficient parallelization, it again comes at the cost of suboptimal cache efficiency: the data loaded for one cell is not revisited until the entire set is processed, potentially leading to increased cache misses.
Therefore maintaining enough parallelism of the smoothing procedure while keeping the method as local as possible are two conflicting goals.   Our algorithm aims to balance between those two aspects: we propose dividing the loop into smaller batches that can be  processed in parallel without simultaneous access to the same data, while the organization of consecutive batches prioritizes data locality. The application requires some tuning as the batch size plays a critical role in terms of overall efficiency.

The remainder of this article is structured as follows. In Section~\ref{sec:multigrid} we briefly describe the multigrid method and its key ingredient, the smoothing step. In Section~\ref{sec:patch-smoothers} we discuss possible implementations of the patch smoothers and in Section~\ref{sec:sequential} we examine a  single use of the smoother. We also explore the impact of different scheduling strategies on the performance of our computations. In Section~\ref{sec:global-operations}, we investigate how the scheduling of tasks, such as grouping neighboring local solvers or applying coloring techniques, can be explited to efficiently exectute the algorithm in thread-parallel enviroment. By analyzing these aspects, we aim to provide insights into the relationship between scheduling choices and performance outcomes.

\section{Geometric finite element multigrid methods}
\label{sec:multigrid}
In this section we describe the mathematical setup of the multigrid method and distinguish between local and global residuals during smoothing. To this end, first
assume a hierarchy of meshes
\begin{gather}
  \mesh_0 \sqsubset \mesh_1 \sqsubset \dots \sqsubset \mesh_L,
\end{gather}
subdividing a domain in $\R^d$, where the symbol
``$\sqsubset$'' indicates nestedness, that is, every cell of mesh
$\mesh_{\ell+1}$ is obtained from a cell of mesh $\mesh_{\ell}$ by refinement. We
note that topological nestedness is sufficient from the algorithmic
point of view, such that domains with curved boundaries can be
covered approximately.

With each mesh $\mesh_\ell$, we associate a finite element space $V_\ell$.
A basis for this space is obtained in the standard way, see for instance~\cite{Ciarlet78}: first, a basis of ``shape functions'' is defined on each cell.
These are associated with ``degrees of freedom'', for instance values in support points for an interpolatory finite element.
If a shape function vanishes on the boundary of the cell, we associate its degree of freedom topologically to the interior of the cell, independent of whether it is a support point or an integral.
If a degree of freedom is topologically on a boundary face/edge/vertex of a cell, it contributes to inter-cell continuity by being identified with degrees of freedom of other cells sharing the same boundary entity.
If the associated shape functions are concatenated, we obtain a basis of $V_\ell$ consisting of continuous functions with local support on cells sharing at least a vertex.

In the following description we will identify $V_\ell$ with the vector space $\mathbb R^{\operatorname{dim}V_\ell}$ of coefficients with respect to this basis. Hence, the Euclidean inner product on $\mathbb R^{\operatorname{dim}V_\ell}$ also serves as an inner product on $V_\ell$.
Equally, we do not distinguish between a finite element function $u_\ell$ and its coefficient vector by notation.
Between these spaces, we introduce transfer operators
\begin{gather}
  \begin{array}{rlcl}
      \rest\ell\colon& V_{\ell+1} &\to&V_\ell,\\
      \prol\ell\colon& V_{\ell} &\to&V_{\ell+1}.
  \end{array}
\end{gather}
As usual, $\prol\ell$ is chosen as the embedding operator of the function spaces and $\rest\ell$ as its
adjoint with respect to the Euclidean inner products in both spaces. Seen as matrices, they are transposes of each other.

We assume that a linear, partial differential equation is defined weakly by a bilinear form $a(\cdot,\cdot)$ and a linear form $b(\cdot)$ for the differential operator and the right hand side.
While the actual equation being solved does not affect the discussion on data structures, our computations are done for Poisson's equation.
The discrete weak formulation of the problem on level $\ell$ is: find $u_\ell\in V_\ell$ such that
\begin{equation}
  \label{eq:weak}
    a (u_\ell,v) =b(v)\quad \forall v\in V_\ell,
\end{equation}
or in matrix notation
\begin{gather}
  \label{eq:matrix}
  A_\ell u_\ell = b_\ell,
\end{gather}
where the matrix $A_\ell$ and the right hand side $b_\ell$ are defined by applying the bilinear form $a(\cdot,\cdot)$ and $b(\cdot)$ to the basis functions of $V_\ell$.

Multigrid methods are efficient solution methods for the discrete linear system \eqref{eq:matrix}, at least for PDE with a regularity gain. When implementing the finite element V-cycle, see for instance~\cite{Bramble93,Hackbusch85}, a so-called smoother is needed on each level in addition to the operators $A_\ell$, $\rest\ell$, and $\prol\ell$. It is usually described by its error propagation operator $\mathcal{S_{\ell}}$.
The simplest example is Richardson's method, where 
  \begin{gather*}
    \mathcal{S_{\ell,\text{Rich}}} = I_\ell - \omega A_\ell.
  \end{gather*}
  Since our goal is a fairly general method not restricted to Poisson's equation, we introduce more powerful smoothers in the next subsection.

\subsection{Vertex patch smoothers}
\label{sec:subspace-correction}

In this article, we focus on vertex patch smoo\-thers \cite{ArnoldFalkWinther00,JanssenKanschat11,KanschatMao15,WitteArndtKanschat21,brubeck2021scalable}.
Since they are based on domain decomposition, they afford sufficient generality to be applied to many different equations, in particular including those involving div and curl operators.
A vertex patch is the subdomain obtained by taking all cells adjacent to a given vertex.
The vertex patch smoother is a subspace correction method with the following structure (we are omitting the level index $\ell$ here, since the smoothers operate on a single level only):
\begin{enumerate}
\item The domain of computation $\Omega$ is covered by vertex patches $\Omega_j$, where $\Omega_j$ is the union of all cells having $x_j$ as a vertex. The indices $j=1,\dots,J$ enumerate all interior vertices of the domain.
Note that this generates an overlapping domain decomposition method.
\item The solution space $V$ corresponding to the discretization on the entire
  mesh is split into subspaces $V_j$ such that the functions in $V_j$ have support in $\Omega_j$. 
  There holds
  \begin{gather}
    \label{eq:main:1}
    V = \sum_{j=1}^J V_j.
  \end{gather}
  With each subspace, we associate a basis $\{v_{j,i}\}_{i=1,\dots,n_j}$, a projection operator $\Pi_j\colon V\to V_j$, and its $\ell_2$-adjoint $\Pi_j^T$.
  \item On each subspace, we employ a local solver, described by the Ritz projection
  \begin{gather}
    \label{eq:main:2}
    P_j v = \Pi_j ^T A_j^{-1} \Pi_j A v.
  \end{gather}
  Here, $A_j$ is the projection of $A$ to the subspace $V_j$.
  Approximate local solvers are covered by the setup albeit the projection property of $P_j$ is lost then. Note that $P_j$ is used to define the error propagation operator, such that it is applied to $v=u^*-u$, where $u^*= A^{-1}b$ is the exact solution and $u$ is its current approximation. Thus, the operation $Av$ actually amounts to computing the residual $b-Au$.
  \item Since additive overlapping subspace correction methods usually perform bad\-ly, we define the smoother as a multiplicative (successive) application of the local solvers represented by the error propagation operator
  \begin{gather}
  \label{eq:s-ell}
      \mathcal{S_{\ell}}=\prod_{j=0}^{J}(I-P_{j}).
  \end{gather}
    Note that the resulting smoother is non-symmetric and depends on the ordering of the subspaces.
\end{enumerate}

In theory, $\Pi_j$ is usually chosen as the $L^2$-projection onto the local space.
In practice, it is the $\ell_2$-projection, namely the selection of the degrees of freedom of $V_j$ from the global space $V$, so that the subspace projection follows the scheme (see also Algorithm~\ref{alg:Loop-naive}) 
\begin{gather}
  \label{eq:main:3}
  \begin{pmatrix}
    u_1\\\vdots\\u_{j_i}\\u_{j_k}\\\vdots\\x_u
  \end{pmatrix}
  \xmapsto{b-Au}
  \begin{pmatrix}
    r_1\\\vdots\\r_{j_i}\\r_{j_k}\\\vdots\\r_n
  \end{pmatrix}
  \xmapsto{\Pi_j}
  \begin{pmatrix}
    r_{j_1}\\\vdots\\r_{{n_j}}
  \end{pmatrix}
  \xRightarrow{A_j^{-1}}
  \begin{pmatrix}
    d_1\\\vdots\\d_{n_j}
  \end{pmatrix}
    \xmapsto{I-\Pi^T_j}
  \begin{pmatrix}
    0\\\vdots\\u_{j_i}\\u_{j_k}\\\vdots \\ 0
  \end{pmatrix}
  .
\end{gather}

\section{Implementation of vertex patch smoothers}
\label{sec:patch-smoothers}

A straight-forward implementation of the successive smoothing operation in equation~\eqref{eq:s-ell} is Algorithm~\ref{alg:Loop-naive}. Given vectors $u$ and $b$ in the finite element space $V$, local corrections are added to $u$. The main ingredient is the function \textsc{\LSolNew}, which projects the residual to the local subspace, solves the local system, and uses the embedding to update the global solution vector $u$.
Since the residual is computed in every step, it has the prohibitive complexity of $\mathcal O(n^2)$ operations and will not be considered any further in this form.
\begin{algorithm}[tp]

\begin{algorithmic}
\State  \Function{\LSolNew}{j,r}
\State{$ r_j = \Pi_j r $} \Comment{Gather}
\State{$ d_j \gets A_j ^{-1} r_j$} \Comment{Solve}
\State \Return{$\Pi_j^T d_j $} \Comment{Scatter}
\EndFunction

\State\For{$j=1,...,N_\text{patches}$}
\State   $ r \gets b-A u$ \Comment{Global residual}
\State   $u \pluseq $ \Call{\LSolNew}{j,r}
\EndFor

\end{algorithmic}

\caption{Subspace correction smoother: straight-forward implementation of the mathematical representation.}
\label{alg:Loop-naive}
\end{algorithm}
In this article, we discuss two ways to restore linear complexity. A common way is ``colorization'', where the index set $1,\dots,J$ is divided into subsets with nonoverlapping patches. This reduces the number of residual computations to a constant independent of the mesh size. It is discussed in subsection~\ref{sec:colorization}. The solution we propose in the present article is a ``local'' computation of residuals and is described first.

\subsection{Local residuals}

Since the local solve in Algorithm~\ref{alg:Loop-naive} only changes a few entries in the result, the residuals of two consecutive steps will also coincide in most entries. Thus, a full matrix vector product in every step is not necessary, as only the changed values of the residual need an update after each local solve. Alternatively, only those entries of the residual $r$ needed for the current local solution must be computed. This algorithm corresponds to the direct implementation of the Gauss--Seidel method.

From the form
\begin{gather*}
    r_{j,i} = f(v_{j,i}) - a(u,v_{j,i}) \qquad i=1,\dots,n_j,
\end{gather*}
we realize that only those degrees of freedom of $u$ coupling with any of the $v_{j,i}$ through $A$ are needed. For conforming finite elements, this is fairly simple: if $V_j$ is the space of functions with support in $\Omega_j$, then the restriction of $u$ to the space $\overline V_j$ consisting of all basis functions whose degrees of freedom are associated with $\Omega_j$ is required. This space is obtained by adding to $V_j$ all basis functions whose node functionals are on the boundary of $\Omega_j$. With this space $\overline V_j$, we associate the embedding operator $\overline\Pi_j^T\colon \overline V_j\to V$ and its $\ell_2$-adjoint projection operator $\overline\Pi_j$.

Thus, we can describe the action of the local solver $P_j$ by the following steps:
\begin{enumerate}
\item \textbf{Gather.} Fetch values of patch-local vector elements from the solution $u$ and right-hand side $b$:
    \begin{gather}
    \label{eq:gather}
    \overline u_j = \overline\Pi_j u\in \overline V_j, \qquad b_j = \Pi_j b \in V_j
  \end{gather}
  This step is generic and corresponds to
  collecting the coefficients of $u_j$ from the global coefficient
  vector of $u$. 
\item \textbf{Evaluate.} Compute local residual using the local operator $\overline A_j$, the restriction of $A$ to $\overline V_j$:
  \begin{gather}
    \label{eq:main:4}
    r_j = b_j - \Pi_j\overline A_j \overline u_j \in V_j.
  \end{gather}
  Here, $\Pi_j$ is implicitly understood as the restriction of the projection $\Pi$ to the subspace $\overline V_j$ of $V$.
\item \textbf{Local solve (LS).} Solve the linear system
  \begin{gather}
    \label{eq:main:5}
    A_j d_j = r_j
  \end{gather}
  to compute the correction $d_j \in V_j$.

\item \textbf{Scatter.} Update the solution $u$:
    \begin{gather}
    \label{eq:scatter}
    u \pluseq \Pi_j^T d_j.
  \end{gather}
  This step is generic and corresponds to
  adding the coefficients $d_{j,i}$ to the respective entries in the global coefficient
  vector of $u$. 
\end{enumerate}

\begin{figure}[tp]
    \centering
    \scriptsize
    \def\svgwidth{\columnwidth}
    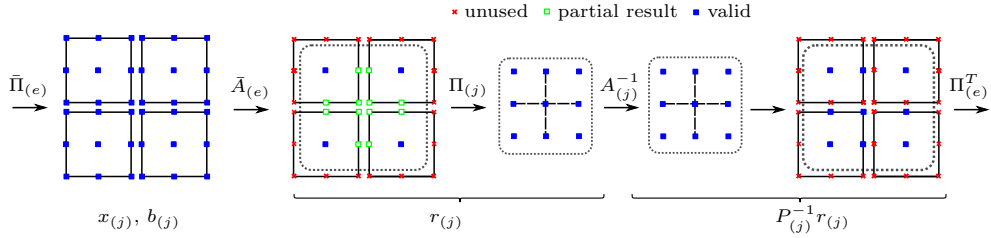
    
    \caption{Steps to evaluate a patch. First, gather data of all involved cells. Second, reinterpret local data for residual computation. Third, restrict to interior of the patch. Fourth, solve in the interior. Fifth, scatter all cell data into the global vector, avoiding double entries.}\label{fig:implementation:steps}
\end{figure}

These steps are visualized in Figure~\ref{fig:implementation:steps}, showing that the algorithm layout is compatible with cell-based matrix-free evaluation~\cite{Kronbichler2012} of $a(u,v_{j,i})$ solely on the cells participating in $\Omega_j$.
Gathering entries from the global data vectors involves all degrees of freedom of all cells of a patch. The action of operator $\overline A$ is computed and subtracted from the local vector cell by cell as well. The values of the residual are rearranged and stored locally to eliminate double values at the interior interfaces and ordered in a way suitable for fast diagonalization. Then, we ignore the boundary values, which contain only partial results without the contributions from cells not in $\Omega_j$, and apply the inverse via fast diagonalization to the interior degrees of freedom on the whole patch (indicated by subscript $j$). Finally, we add the local solution to the global data vector, where we make sure we only update relevant values.

The version of Algorithm~\ref{alg:Loop-naive} with local residuals can be found in Algorithm~\ref{alg:Loop-sequential}. The four steps described above are encapsulated in the function \textsc{\LUpNew}. We refer to this method as the \emph{combined} algorithm as local residual computation and application of the local solvers are combined in a single loop.
\begin{algorithm}[tp]

\begin{algorithmic}

  \State \Function{\LUpNew}{j,u}
  \State    $ u_j \gets \overline\Pi_j u$ \Comment*{Gather}
  \State    $ r_{j} \gets \Pi_j b - \Pi_j\overline A_j u_j$ \Comment*{Evaluate}
  \State    $ d_j \gets A_i ^{-1} r_j$ \Comment*{Solve}
  \State \Return{ $\Pi_j^T d_j$ }\Comment*{Scatter}
  \EndFunction
  
\State \For{$j=1,\ldots,N_\text{patches}$}
     \State $u \pluseq$ \Call{\LUpNew}{$j$,$u$}
    
    \EndFor
  
\end{algorithmic}
  \caption{Combined smoother: a single patch loop where local residuals are combined with local solvers.}
  \label{alg:Loop-sequential}
\end{algorithm}

\subsection{Colorization}
\label{sec:colorization}

Algorithm~\ref{alg:Loop-naive} requires the application of the operator $A$ for every local solution step of the multiplicative algorithm~\eqref{eq:s-ell}, which is highly inefficient.
It is known that additive domain decomposition methods only need a single application of $A$, see also~\cite{munch2023cache}.
Unfortunately, additive overlapping methods involve a considerably higher number of iterations, see for instance~\cite{WitteArndtKanschat21}.
This problem can be mitigated by ``colorization'', which goes back to the red-black coloring of the Gauss--Seidel algorithm.

Each local update reads the data from the whole patch including the boundary, but only writes data corresponding to the interior degrees of freedom.
Therefore, a set of mutually non-overlapping vertex patches can be operated on in any order without changing the result.
This leads to the following form of the smoother:
\begin{gather}
  \label{eq:colorized}
  \mathcal{S_{\ell}}
  =\prod_{k}^{N_c}\left(I- \left(
    \sum_{j\in C_{k}}\Pi_{j}^{T}A_{j}^{-1}\Pi_{j}
    \right)A\right),
\end{gather}
where the index sets $C_k$, $k=1,...,\ncol$, are called \emph{colors}. If they form sets of subspaces which do not overlap, that is,
\begin{gather}
    \Omega_i\cap\Omega_j = \emptyset \qquad i\neq j, \quad i,j\in C_k,
\end{gather}
the two versions~\eqref{eq:s-ell} and~\eqref{eq:colorized} are mathematically equivalent. Algorithm~\ref{alg:Loop_by_colors} shows an implementation of~\eqref{eq:colorized}.
\begin{algorithm}[tp]

\begin{algorithmic}
  \For{$k=1,...,\ncol$
  }
   \State $ r \gets b-A u$
    
   \State \ForAll{$ j\in C_k$
    }
    \State  $u \pluseq $ \Call{\LSolNew}{j,r}
    \EndFor
  
\EndFor
\end{algorithmic}
\caption{Separated, colorized smoother}
\label{alg:Loop_by_colors}
\end{algorithm}
It requires computing the action of the global operator $A$ only $\ncol$ times, which is separated from the loop over the patch solvers. Hence, we will refer to this option as the \emph{separated, colorized} algorithm.
In Algorithm~\ref{alg:Combined_by_colors} below, we combine local computation of the residuals with colorization.
\begin{algorithm}[tp]

\begin{algorithmic}
  \For{$k=1,...,\ncol$}
    \ForAll{$ j\in C_k$}
      $u \pluseq $ \Call{\LUpNew}{j,u}
    \EndFor
    \EndFor
  
\end{algorithmic}
\caption{Combined, colorized smoother.}
\label{alg:Combined_by_colors}
\end{algorithm}

\subsection{Comparison of sequential implementations}
\label{sec:sequential}

We have now three implementations of the same mathematical algorithm, which we expect to be equivalent up to the order in which vertex patches are traversed. We now show that, while employing a very similar number of floating point operations, they have quite different behavior when it comes to memory access. The latter in particular is important for efficient implementation on modern hardware.
We note that Algorithm~\ref{alg:Loop-naive} with global matrix-vector products and without colorization is inherently inefficient and not considered here.

On a machine equipped with two AMD EPYC 7282 processors (16 cores, 2.8 GHz, 64 MBytes L3 cache, using SIMD of the AVX2 instruction set), we run the algorithms on a single core, which leads to a core-bound scenario. This has been verified by measuring the run-time while omitting the majority of computations, which leads to a significantly reduced run time. Hence, we may use the sequential run time as a relative measure of the number of arithmetic operations and indirect addressing costs.

When we say that the execution of the algorithms is ``core-bound'', we describe a scenario where the performance is limited by in-core performance, which for our code is a mix of saturation of compute units in the sum factorization part, saturation of L1/L2 cache bandwidth for the operation at quadrature points, saturation of load/store units with respective buffers for the indirect addressing in the gather/scatter steps, or pipeline stalls due to latency of dependent operations. In parallel, when many cores share the memory interface, the relative cost of memory transfer increases and the bandwidth of memory controllers to the main memory can become the limiting aspect instead. The latter term is abbreviated as ``memory-bound'' in the following.

\begin{figure}[tp]
\begin{centering}
    \includegraphics[width=.7\textwidth]{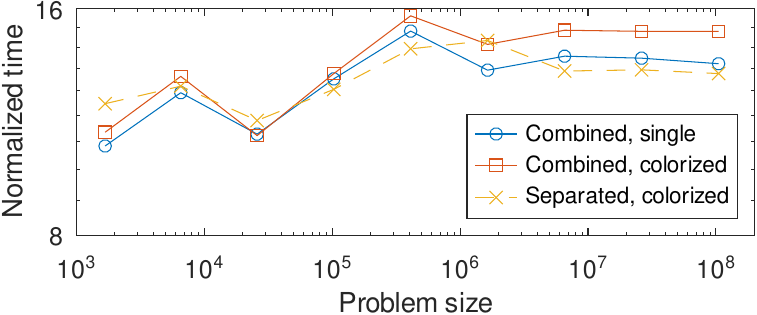}
\caption{Dependence of average time for one smoothing step on mesh refinement with the
combined (Algorithm~\ref{alg:Loop-sequential}), combined colorized  (Algorithm~\ref{alg:Combined_by_colors}) and separated colorized (Algorithm~\ref{alg:Loop_by_colors}) strategies. $\mathbb{Q}_5$ shape functions in two dimensions with Z-curve ordering. The timings are normalized by the time for one optimized matrix-vector multiplication.
}
\label{fig:serial-mesh}
\end{centering}
\end{figure}

We anticipate that the computational expense of the operation in a two-dimen\-sio\-nal~(2D) problem should be at least eight times higher than that required for computing the matrix-vector product. This estimation is based on the fact that each cell is visited four times (8 times in 3D), where two operations are performed each time: the matrix-vector product and the local solver. The latter involves somewhat more arithmetic operations due to the coupling of all the unknowns within the patch within sum factorization. In Figure~\ref{fig:serial-mesh}, we show the measured run time for one smoothing step with each algorithm depending on the mesh size, normalized by the time to compute the matrix-vector product with matrix-free operator evaluation \cite{Kronbichler2012}. The result verifies the performance expectations with a cost of eight to sixteen matrix-vector products.
All curves are more or less flat, indicating the expected linear dependence on the number of degrees of freedom. Furthermore, the difference between the different implementations does not exceed 10\%, indicating that the computational effort is almost the same. The combined versions show a slightly higher compute time because the residual around shared degrees of freedom involves additional arithmetic operations and the data access of cells is vectorized more efficiently for the global matrix-vector product~\cite{Kronbichler2012}. 
For the two combined versions, this is expected. For the separated version, we point out that these experiments are run on uniform meshes. Thus, the coloring can be done in an almost optimal way, such that the patches in each color cover almost the complete domain. Hence, almost all degrees of freedom have been modified when the next global residual is computed.
We conclude that in our experimental settings, the computational effort of the different strategies of the patch smoothing algorithm is almost the same.

Next we investigate the memory transfer between main memory and the L3 cache.
\begin{figure}
\begin{centering}
    \includegraphics[width=.48\textwidth]{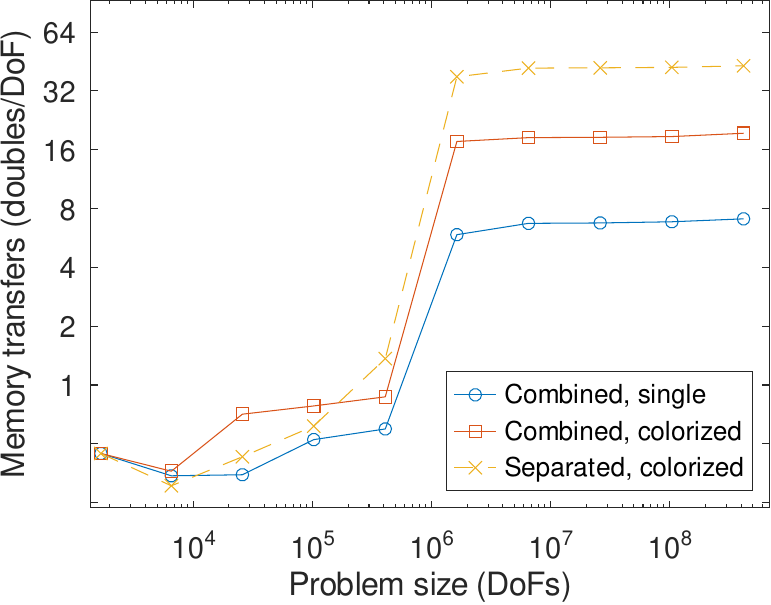}
    \hfill
        \includegraphics[width=.48\textwidth]{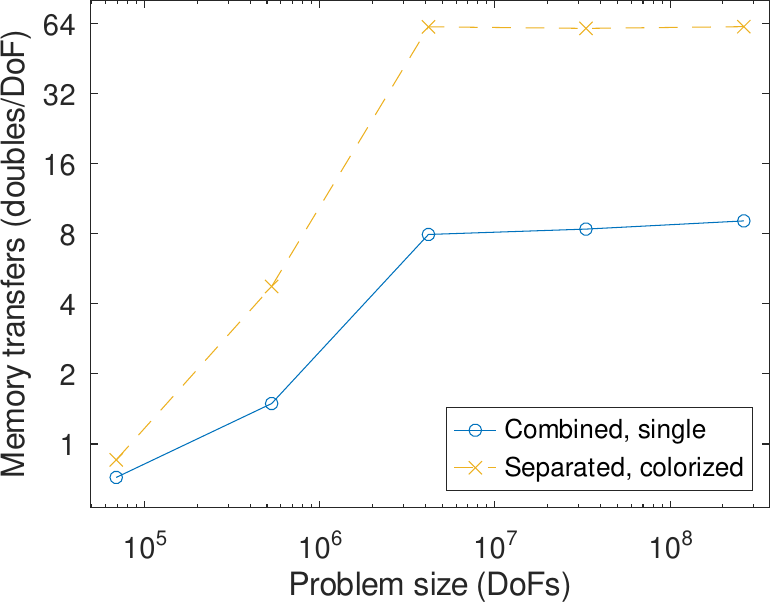}
\caption{Measured transfer between main memory and L3 cache in each smoothing step averaged over 100 steps. 2D on the left, 3D on the right, $\mathbb Q_5$ elements.
}
\label{fig:serial-mesh-liqwid}
\end{centering}
\end{figure}
Figure~\ref{fig:serial-mesh-liqwid} shows the measured transfer for a single smoothing step as a function of the mesh size, displayed as the number of double variables (8 bytes) per DoF. We use the same machine and run the algorithm for biquintic elements on two-dimensional uniform meshes.
We see that the L3 cache size is exhausted around one million degrees of freedom. Above this threshold, the amount of memory transfers per degree of freedom, measuring the cache efficiency of the algorithm, is almost constant for all versions.

Hence, we compare the values in the rightmost data points of Figure~\ref{fig:serial-mesh-liqwid} in more detail in Table~\ref{tab:L3-cache}. Note that the problem size exceeds the cache size by two orders of magnitude, such that the numbers can be considered representative.
\begin{table}[tp]
  \centering
  \caption{Transfer between main memory and L3 cache for Algorithms 2--4 represented as the number \texttt{double} variables per degree of freedom. The mesh size is such that problem sizes are about $10^8$ degrees of freedom.}
  \label{tab:L3-cache}
  \begin{tabular}{llc|ccc}
    \toprule
    & FE &ordering  & comb. single & comb. col. & sep. col. \\
    \cmidrule[\lightrulewidth]{1-6}
    \multirow{4}{*}{2D}
&  \multirow{2}{*}{$\mathbb Q_5$}
& Z-curve &  7.1 &19.4  & 43.0 \\
&   & hierarchical  & 32.9 & 38.4 &  61.5\\
    \cmidrule[\lightrulewidth]{2-6}
  &   \multirow{2}{*}{  $\mathbb Q_3$ }
  & Z-curve &8.1     & 21.9 &   46.8  \\
  &    & hierarchical & 47.0&   51.3 &     74.8  \\
    \midrule[\lightrulewidth]
    \multirow{1}{*}{3D}
    &   \multirow{2}{*}{ $\mathbb Q_5$}
& Z-curve &9.1 &   47.2 &  62.2 \\
&  & hierarchical  & 44.4 & 80.0 &  94.6  
     \\
    \bottomrule
  \end{tabular}
\end{table}
Here, another issue arises: While ordering of the patches does not matter for assessing the number of operations, it becomes important when we study the cache efficiency. Therefore, we include the effect of the order of the vertex traversal when comparing the different algorithms discussed above.
We consider two vertex patch orderings. Ordering vertices along a Z-curve traversing the mesh cells of the current level affords fairly good locality~\cite{bangerth2012algorithms,kronbichler2022cg}. Particularly bad data locality is obtained by following the hierarchical construction of the same Z-curve, where first the vertices on the coarsest mesh level are numbered, then those on the next and so on, resulting in a high number of long distance jumps.

We see that the combined algorithm with a single patch loop (Algorithm~\ref{alg:Loop-sequential}) and Z-curve ordering has by far the lowest memory transfer.
The measured memory transfer is still higher than the minimum of three doubles per unknown (read $b$, read and write $u$), because the indices of the unknowns for both the localized matrix-vector product and for the patch gather and scatter operators $\Pi_j$ also need to be accessed, as well as some metadata for the application of $A_j^{-1}$. Nevertheless, the overhead is about a factor of two.

This is also the algorithm most susceptible to deterioration due to inefficient ordering, as the memory transfer increases by factors between 4 and 5 for hierarchical order.  Additionally, it is also an inherently sequential algorithm.

Thus, we turn our attention to the second column, the combined colorized algorithm and analyze the price to be paid for conflict-free read and write operations. We first discuss the Z-curve ordering. Since the loop over colors traverses the whole vector four times in two dimensions, we also expect this factor in the memory transfers. The observed factor is slightly less at about three, as in fact not all data can be cached in the combined algorithm with single patch loop. In three dimensions, the transfer ratio is close to twice as high, corresponding to the fact that we have eight colors instead of four. Since the memory layout is not as efficient as in the single loop, the effect of changing the ordering from Z-curve to hierarchical is less pronounced here, but still indicates that a reasonably good ordering is essential.

The separate implementation of smoothing and residual computation, namely Algorithm~\ref{alg:Loop_by_colors}, is provided in the third column for reference. As expected, memory transfer rates are considerably higher than with the combined algorithm.

\section{Parallel application of the smoother}
\label{sec:global-operations}

One of the main concerns for parallel implementation is the presence of data races, which occur when several threads access and modify a shared data structure simultaneously.
In the context of patch smoothing, data races can only occur when multiple threads work on the same regions of the vector $u$ by acting on overlapping patches simultaneously.
The additive part of the colorized smoother~\eqref{eq:colorized} can be executed in parallel without causing race conditions and thus does not require synchronization.
This comes at the cost of cache efficiency since the cache capacity is typically exhausted while traversing one color in its entirety, implying and additional load operation for the next color.
Our goal here is to develop an algorithm which is more cache efficient and free of race conditions.

\subsection{Batch processing}

In order to maintain data locality and prevent potential race conditions in the proposed multigrid patch smoothers, the patches of the mesh are enumerated and colored in a specific way. The patches, denoted as $[1, \dots, N_\text{patches}]$, are first ordered in such a way that data access during loops is optimized, that is, following a Z-curve. Then, the patches are grouped into colors to prevent data race conditions, while preserving the global order within each color. Finally, the  colored patches are further divided into batches of fixed size that are supposed to fit into the L3 cache. To summarize, the procedure used to set up the loop structure, depending on the parameter $n_B$ that allow us to tune the batch size, is done in the following steps:
\begin{enumerate}
\item Generate global patch ordering following a Z-curve.
\item Group patches by colors, but keep the ordering within each color.
\item Divide each color $k$ into batches $P_{b, k}$, $b= 1,..., N_{\text{batches}}$, of prescribed size $n_B$. The last batch in each color might contain a lower number of patches.
\end{enumerate}

 The resulting batches can then be processed in the smoother application, as demonstrated in Algorithm~\ref{alg:Loop-parallel}. The outermost loop is over batches, the middle loop is over colors, and the innermost parallel loop, which is implemented using \texttt{TBB} parallel \texttt{for} loops, iterates over patches. This approach results in an execution that is equivalent to three nested loops, and allows for efficient parallel processing of the patches. Note that there is a synchronization point after each color in the second loop, note also the discussion on synchronization in conjunction with coloring in \cite[Sec.~III]{Kormann2011}.

\begin{algorithm}[tp]
\begin{algorithmic}
  \For{$b= 1,\dots,N_{\text{batches}}$ 
  }
    \For{$k=1,...,\ncol$
    }
      \ForAll{$ j\in P_{b,k}$ 
      }
        \Call{\LUpNew}{j,u}
      \EndFor
    \EndFor
  \EndFor
\end{algorithmic}
\caption{Batch processing, combined, colorized method}
\label{alg:Loop-parallel}
\end{algorithm}

The size of each batch plays a crucial role in the performance of the proposed multigrid patch smoothers. A large number of small batches implies a high number of synchronization points between threads. This cost as well as a small number of patches within $P_{b,k}$ limit the available parallelism. On the other hand, if batches are too large, the data accessed from one color to the next will exceed the cache capacity and thus increase memory transfer. In particular, the extreme case where batches consist of a single patch results in sequential execution, while batches corresponding to whole colors are the most suitable for parallel processing but with a decreased data locality. The optimal size of batches aims to find a sweet spot to balance between these two extremes, and should be a multiple of the number of available threads to avoid workload imbalances.
It is worth noting that the batches can be considered as new `colors' and the loop is the same as in Algorithm~\ref{alg:Loop-parallel}.

By partitioning the patches in a way that their data fits into the L3 cache and assigning work from the same partition to threads, we can optimize the data access by reducing the number of times data needs to be loaded from (slow) main memory, ideally once per NUMA domain. This partitioning reduces the number of cells to be processed and touched in parallel across threads when passing from one color to next (where there is a considerable number of shared data), making the working set fit into caches.

\begin{figure}
    \centering
    \includegraphics[width=.6\textwidth]{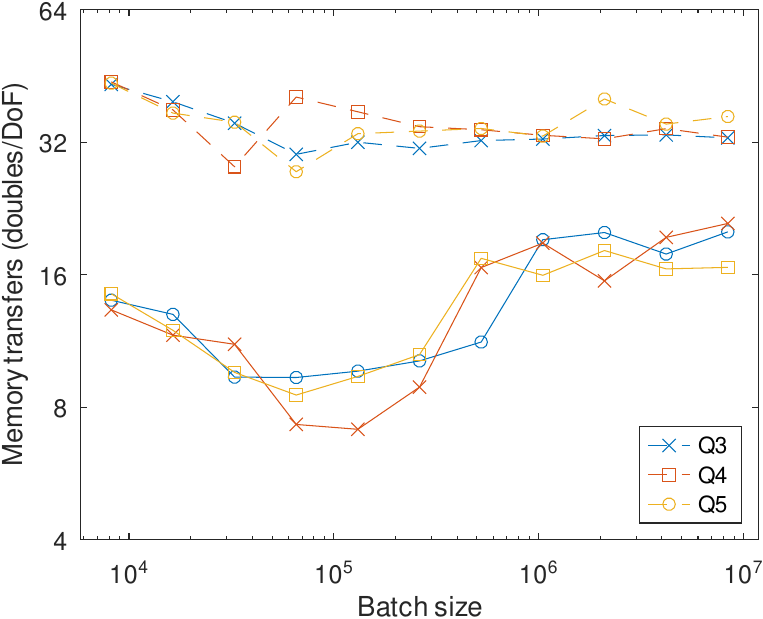}
    \caption{Thread-parallel application with 64 threads: Memory transfer per degree of freedom vs. batch size (expressed as the number of patches). Timings with the separated algorithms are shown with dashed lines, the combined application with solid lines.  Fixed mesh with 10 refinements in 2D.}
                                
    \label{fig:theards_mem_transfer}
\end{figure}

\begin{figure}
    \centering
    \includegraphics[width=.6\textwidth]{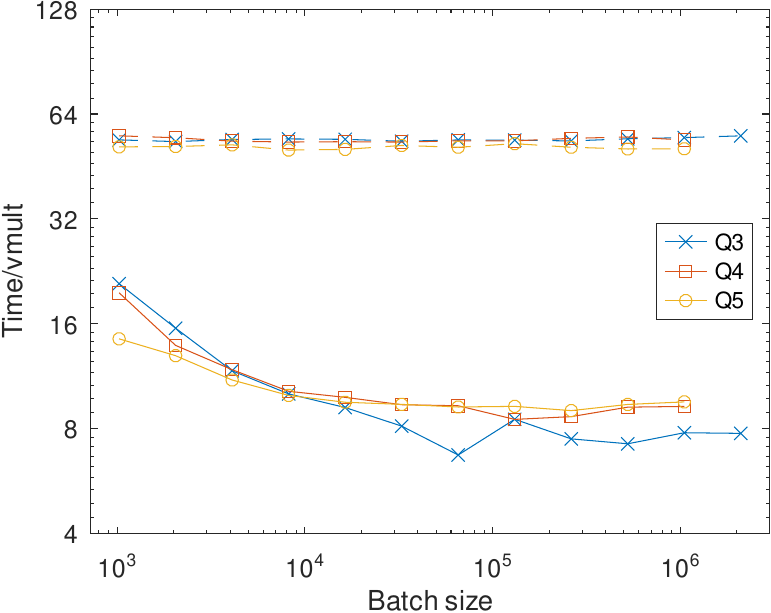}
    \caption{Thread-parallel application  with 32 cores Time/DoF vs batch size. Fixed mesh 2D with 10 refinements, corresponding to the rightmost point from Figure~\ref{fig:serial-mesh-liqwid}.
    Solid is combined, dashed is separated.
    }
    \label{fig:theard_perf_task}
\end{figure}

\paragraph{Comparison of parallel implementations}
In order to evaluate the performance of our parallel patch-smoother implementation, we measure both the execution time and memory transfers, similar to our evaluation of the sequential implementation. However, unlike the sequential implementation, we use a fixed mesh with 10 refinements in two dimensions to evaluate how the batch size affects the performance.
The results are shown in Figure~\ref{fig:theards_mem_transfer} and Figure~\ref{fig:theard_perf_task}.

Figure~\ref{fig:theards_mem_transfer} demonstrates the average memory transfers per smoothing step versus batch size, the latter ranging from  8k to 8M in powers of 2. The transfers are lower for the proposed algorithm (for all batch sizes). There is a pronounced range where transfers are as lowest, which is between $32768$ and $262144$. Compared to very large batches, transfers are around 3 times lower in this region. The 2-fold difference in memory transfer between the combined and separated applications for the largest batch size (corresponding to \textit{color-by-color} execution) can be attributed to the computation of the residual being combined with the application of the inverse. This results in all operations being local, reducing the memory transfer.

Figure~\ref{fig:theard_perf_task} shows the average execution time versus batch size,
normalized against the time of a single matrix-vector product. 
As the batch size increases, the execution time decreases,
indicating a reduction in the cost of thread synchronization.
Furthermore, we observe a shift towards lower memory transfer.
To confirm this conjecture, we performed a similar experiment omitting most computations on each patch, hence measuring the time for memory access; we observed greatly reduced timing for the batched version.
However, the effect of batch size on performance is not very significant as the method has already become core-bound due to the combination of residual computation and inverse application being performed together. For this reason, the synchronization of threads between processing batches still plays a critical role in determining the final execution time.  In high-end machines equipped with CPUs featuring a less generous memory bandwidth relative to their computational capabilities, we anticipate the would be similar to the one depicted in Figure~\ref{fig:theards_mem_transfer}. 

\begin{figure}
    \centering
    \includegraphics[width=.3\textwidth]{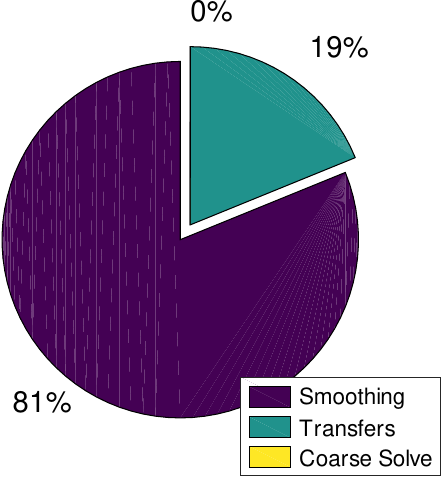}
    \caption{Thread-parallel application in 3D: time consumed by subprocedures of multigrid method. Fixed mesh with 5 refinements Q5
        }
    \label{fig:mg_ingredients}
\end{figure}

Finally, we present the performance analysis of our multigrid solver in Table~\ref{tab:threads-perf}. To provide a comprehensive view of the solver's behavior, we depict the time distribution across various subprocedures using the pie chart shown in Figure~\ref{fig:mg_ingredients}. 
The smoothing still dominates the computational workload, while the accounts for a significant 19\% of the total execution time.  The coarse solve step consumes an almost imperceptible amount of time in the context of the entire procedure.

\begin{figure}
    \centering
    \includegraphics[width=.6\textwidth]{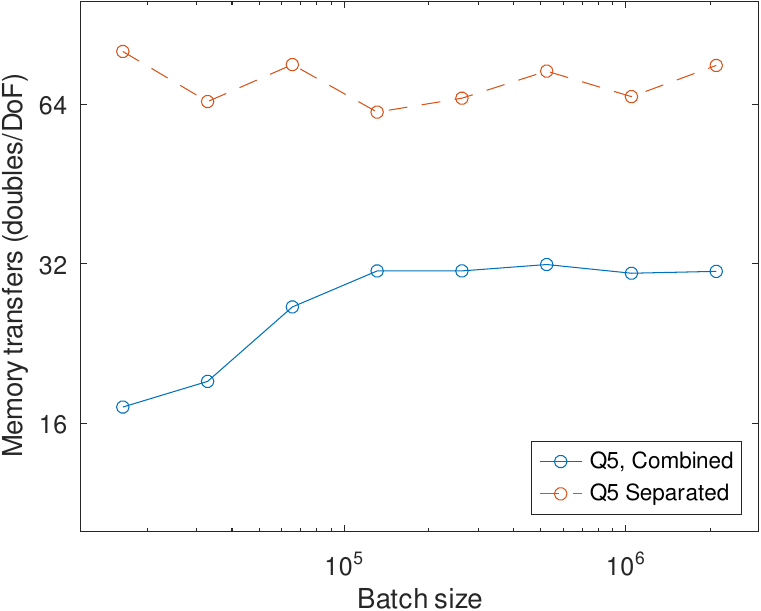}
    \caption{Thread-parallel application in 3D: Memory transfer per DoF vs number of patches in batch. Fixed mesh with 5 refinements Q5.}
    \label{fig:theards_mem_transfer_3D}
\end{figure}

\begin{table}[tp]
  \centering
  \caption{Best average timing for Algorithms~\ref{alg:Loop-parallel} on 32 cores. Residual reduction by $10^{-12}$. Mesh: single patch refined 10 times in 2D and 6 in 3D
    \label{tab:threads-perf}}
\begin{tabular}{lccc}
\hline 
&\multicolumn{2}{c}{2D} & 3D\tabularnewline

&$\mathbb{Q}_3$ & $\mathbb{Q}_5$ & $\mathbb{Q}_5$ \tabularnewline
\hline 
Problem size  &37,761,025 & 104,878,081 & 263,374,721\tabularnewline
Smoothing  &0.197856s & 0.387558s & 2.37646s\tabularnewline
vmult  &0.034349s & 0.0812562s & 0.213686s\tabularnewline
Iterations &6 & 5 & 6\tabularnewline
Solving &12.0369s & 23.3042s & 115.17s\tabularnewline

\hline 
\end{tabular}
\end{table}

\section{Conclusion}

In conclusion, we have focused on improving the performance of patch smoothers, which have been shown to be memory-bound. By reducing memory transfers, we have improved the performance of these smoothers, making them a more efficient choice for multigrid methods.

We have achieved this by combining the computation of the residual with the application of the smoother making the smoother perform similarly to a Gauss--Seidel method.  Our results show that this approach leads to a two-fold reduction in memory transfer. Additionally, by properly ordering the patches in a sequential implementation, we have been able to further reduce memory transfers.

In our thread-parallel implementation, the sequence of patches is divided into batches for parallel processing. 
Therefore, it is important to balance the batch size to achieve the optimal trade-off between minimizing the number of synchronizations and minimizing memory transfer. In our experiments, conducted using a mid-range processor, we found that the algorithm is compute-bound for all batch sizes, but we expect that CPU computing capabilities will continue to grow faster than memory bandwidth in the future, changing the balance.

\FloatBarrier

\appendix

\section{Implementation and test problem}

\subsection{Evaluation of the combined operator}

We have implemented the algorithm using the \texttt{deal.II} finite element
library~\cite{dealII94,dealii2019design} and its matrix-free infrastructure~\cite{Kronbichler2012, Kronbichler2017a}. 
This infrastructure loops 
over cells and processes multiple cells (cell batches) at once in a
vectorization-over-cells manner. In our case, we loop manually over the patches
in a sequence defined by Algorithm~\ref{alg:Loop-parallel}. 
To be able, nevertheless, to use as much infrastructure of base finite element
library and to vectorize operations, we
compute residuals of cells within a patch at once, by using \texttt{deal.II}'s
capabilities to define cell batches on the fly~\cite{dealII94}. We obtain the following
extended algorithm: 
\begin{enumerate}
\item gather DoFs of each cell of a patch,
\item evaluate the residual contributions of each cell,
\item sum up the residual contributions and store the data in lexicographical order,
\item apply the local smoother,
\item split up the result according to the cells, and
\item scatter the results back into the global vector, using a mask to avoid 
adding contributions multiple times.
\end{enumerate}
One could merge the last two steps, however, we decided to add the intermediate
step (5) to be able to reuse for scattering infrastructure used during gathering (1) 
and to not load additional indices.
The steps are visualized in Figure~\ref{fig:implementation:steps}. It shows, particularly, valid and
invalid vector entries, which arise due to the fact that only residuals in the 
interior of the patches can be computed and due to the masking.

To evaluate the residual, we use numerical quadrature, which involves the basis
change from support points to quadrature points and the computation of derived
quantities, which we perform via sum factorization~\cite{Melenk99fullydiscrete, Orszag1980}. Furthermore, we need
certain mapping information at the quadrature points, e.g., for a Laplace operator
the inverse Jacobian and its determinant. Since we are concentrating on affine meshes
in this publication, the mapping information at each quadrature point is the
same and can be easily cached. All this infrastructure is provided 
by \texttt{deal.II}~\cite{Kronbichler2012}.

\subsection{Test problem}
We test our smoother on the Poisson equation
\begin{eqnarray*}
\Delta u & =f & \text{in \ensuremath{\Omega}}\\
u & =0 & \text{on \ensuremath{\partial\Omega} }
\end{eqnarray*}
on a unit hypercube domain  $ [0,1]^d=\Omega\subset\mathbb{R}^{d}$ , where $d=2,3$. We consider a uniform, Cartesian mesh so that we obtain a problem with a tensor-product structure, that could be exploited to provide efficient local solvers. 

\subsection{Local inverse}

In our experiments, we apply the fast diagonlization method~\cite{Lynch1964} to compute
the inverse of the system matrix restricted to the interior of the vertex-patch or 
its approximation. The motivation for fast diagonalization methods is that,
for Cartesian meshes, the matrix of the $d$-dimensional Laplace operator 
can be expressed as the
tensor product of one dimensional mass matrices $M_i$ and stiffness matrices $K_i$ (here shown for 3D):
\begin{align*}
A_b &= M_2 \otimes M_1\otimes K_0 + M_2 \otimes K_1\otimes M_0 + K_2 \otimes M_1\otimes M_0
\end{align*}
This can be transformed to:
\begin{align}\label{eq:fdm:eigen}
A_b  &= T_2 \otimes T_1 \otimes T_0 (\Lambda_2 \otimes I \otimes I + I \otimes \Lambda_1 \otimes I + I \otimes I \otimes \Lambda_0 ) T^T_2 \otimes T^T_1 \otimes T^T_0,
\end{align}
with $T_i$ and $\Lambda_i$ being the eigenvectors and eigenvalues obtained from
a generalized eigendecomposition of $K_i T_i = \Lambda_i M_i T_i$. 
The inverse of \eqref{eq:fdm:eigen} is explicitly given by
\begin{align}\label{eq:fdm:inverse}
A_b^{-1}=T_2 \otimes T_1 \otimes T_0 (\Lambda_2 \otimes I \otimes I + I \otimes \Lambda_1 \otimes I + I \otimes I \otimes \Lambda_0 )^{-1} T^T_2 \otimes T^T_1 \otimes T^T_0.
\end{align}
The evaluation of~\eqref{eq:fdm:inverse} can be implemented computationally
efficiently again by using sum factorization. Fast diagonalization methods have
been also applied, e.g., by~\cite{brubeck2021scalable,munch2023cache,WitteArndtKanschat21} as domain
solver in the context of vertex-patch smoothers.

\FloatBarrier


%% file: vmult_vs_new.pdf_tex
\begingroup%
  \makeatletter%
  \providecommand\color[2][]{%
    \errmessage{(Inkscape) Color is used for the text in Inkscape, but the package 'color.sty' is not loaded}%
    \renewcommand\color[2][]{}%
  }%
  \providecommand\transparent[1]{%
    \errmessage{(Inkscape) Transparency is used (non-zero) for the text in Inkscape, but the package 'transparent.sty' is not loaded}%
    \renewcommand\transparent[1]{}%
  }%
  \providecommand\rotatebox[2]{#2}%
  \newcommand*\fsize{\dimexpr\f@size pt\relax}%
  \newcommand*\lineheight[1]{\fontsize{\fsize}{#1\fsize}\selectfont}%
  \ifx\svgwidth\undefined%
    \setlength{\unitlength}{1824.93480317bp}%
    \ifx\svgscale\undefined%
      \relax%
    \else%
      \setlength{\unitlength}{\unitlength * \real{\svgscale}}%
    \fi%
  \else%
    \setlength{\unitlength}{\svgwidth}%
  \fi%
  \global\let\svgwidth\undefined%
  \global\let\svgscale\undefined%
  \makeatother%
  \begin{picture}(1,0.21967449)%
    \lineheight{1}%
    \setlength\tabcolsep{0pt}%
    \put(0,0){\includegraphics[width=\unitlength,page=1]{vmult_vs_new.pdf}}%
    \put(0.4457509,0.00103385){\color[rgb]{0,0,0}\makebox(0,0)[t]{\lineheight{0}\smash{\begin{tabular}[t]{c}$r_{(j)}$\end{tabular}}}}%
    \put(0,0){\includegraphics[width=\unitlength,page=2]{vmult_vs_new.pdf}}%
    \put(0.81973686,0.00103385){\color[rgb]{0,0,0}\makebox(0,0)[t]{\lineheight{0}\smash{\begin{tabular}[t]{c}$P^{-1}_{(j)} r_{(j)}$\end{tabular}}}}%
    \put(0,0){\includegraphics[width=\unitlength,page=3]{vmult_vs_new.pdf}}%
    \put(0.97973616,0.13578101){\color[rgb]{0,0,0}\makebox(0,0)[t]{\lineheight{0}\smash{\begin{tabular}[t]{c}$\Pi_{(e)}^T$\end{tabular}}}}%
    \put(0,0){\includegraphics[width=\unitlength,page=4]{vmult_vs_new.pdf}}%
    \put(0.62275209,0.13578101){\color[rgb]{0,0,0}\makebox(0,0)[t]{\lineheight{0}\smash{\begin{tabular}[t]{c}$A^{-1}_{(j)}$\end{tabular}}}}%
    \put(0.46793632,0.13578101){\color[rgb]{0,0,0}\makebox(0,0)[t]{\lineheight{0}\smash{\begin{tabular}[t]{c}$\Pi_{(j)}$\end{tabular}}}}%
    \put(0,0){\includegraphics[width=\unitlength,page=5]{vmult_vs_new.pdf}}%
    \put(0.01915307,0.13578101){\color[rgb]{0,0,0}\makebox(0,0)[t]{\lineheight{0}\smash{\begin{tabular}[t]{c}$\bar{\Pi}_{(e)}$\end{tabular}}}}%
    \put(0.13186556,0.00103385){\color[rgb]{0,0,0}\makebox(0,0)[t]{\lineheight{0}\smash{\begin{tabular}[t]{c}$x_{(j)}$, $b_{(j)}$\end{tabular}}}}%
    \put(0,0){\includegraphics[width=\unitlength,page=6]{vmult_vs_new.pdf}}%
    \put(0.46411137,0.21120293){\color[rgb]{0,0,0}\makebox(0,0)[lt]{\lineheight{0}\smash{\begin{tabular}[t]{l}unused\end{tabular}}}}%
    \put(0.55954723,0.21120293){\color[rgb]{0,0,0}\makebox(0,0)[lt]{\lineheight{0}\smash{\begin{tabular}[t]{l}partial result\end{tabular}}}}%
    \put(0.71233692,0.21120293){\color[rgb]{0,0,0}\makebox(0,0)[lt]{\lineheight{0}\smash{\begin{tabular}[t]{l}valid\end{tabular}}}}%
    \put(0,0){\includegraphics[width=\unitlength,page=7]{vmult_vs_new.pdf}}%
    \put(0.24535992,0.13395018){\color[rgb]{0,0,0}\makebox(0,0)[t]{\lineheight{0}\smash{\begin{tabular}[t]{c}$\bar{A}_{(e)}$\end{tabular}}}}%
  \end{picture}%
\endgroup%

%% file: main.bbl
\begin{thebibliography}{10}

\bibitem{dealII94}
{\sc D.~Arndt, W.~Bangerth, M.~Bergbauer, M.~Feder, M.~Fehling, J.~Heinz,
  T.~Heister, L.~Heltai, M.~Kronbichler, M.~Maier, P.~Munch, J.-P. Pelteret,
  B.~Turcksin, D.~Wells, and S.~Zampini}, {\em The \texttt{deal.II} library,
  version 9.5}, J. Numer. Math., 31 (2023), pp.~231--246,
  \url{https://doi.org/10.1515/jnma-2023-0089},
  \url{https://dealii.org/deal95-preprint.pdf}.

\bibitem{dealii2019design}
{\sc D.~Arndt, W.~Bangerth, D.~Davydov, T.~Heister, L.~Heltai, M.~Kronbichler,
  M.~Maier, J.-P. Pelteret, B.~Turcksin, and D.~Wells}, {\em The {deal.II}
  finite element library: Design, features, and insights}, Comput. \& Math.
  Appl., 81 (2021), pp.~407--422,
  \url{https://doi.org/10.1016/j.camwa.2020.02.022},
  \url{https://arxiv.org/abs/1910.13247}.

\bibitem{ArnoldFalkWinther97Hdiv}
{\sc D.~N. Arnold, R.~S. Falk, and R.~Winther}, {\em Preconditioning in
  {$H({\rm div})$} and applications}, Math. Comput., 66 (1997), pp.~957--984,
  \url{https://doi.org/10.1090/S0025-5718-97-00826-0}.

\bibitem{ArnoldFalkWinther00}
{\sc D.~N. Arnold, R.~S. Falk, and R.~Winther}, {\em Multigrid in {$H({\rm
  div})$} and {$H({\rm curl})$}}, Numer. Math., 85 (2000), pp.~197--217,
  \url{https://doi.org/10.1007/PL00005386}.

\bibitem{bangerth2012algorithms}
{\sc W.~Bangerth, C.~Burstedde, T.~Heister, and M.~Kronbichler}, {\em
  Algorithms and data structures for massively parallel generic adaptive finite
  element codes}, ACM Trans. Math. Softw., 38 (2012), pp.~14/1--28,
  \url{https://doi.org/10.1145/2049673.2049678}.

\bibitem{Bramble93}
{\sc J.~H. Bramble}, {\em Multigrid Methods}, no.~294 in Pitman research notes
  in mathematics series, Longman Scientific, 1993.

\bibitem{brubeck2021scalable}
{\sc P.~D. Brubeck and P.~E. Farrell}, {\em A scalable and robust vertex-star
  relaxation for high-order {FEM}}, SIAM J. Sci. Comput., 44 (2022),
  pp.~A2991--A3017, \url{https://doi.org/10.1137/21M1444187}.

\bibitem{Ciarlet78}
{\sc P.~G. Ciarlet}, {\em The Finite Element Method for Elliptic Problems},
  North-Holland, 1978.

\bibitem{Dumbser2018}
{\sc M.~Dumbser, F.~Fambri, M.~Tavelli, M.~Bader, and T.~Weinzierl}, {\em
  Efficient implementation of {ADER} discontinuous {G}alerkin schemes for a
  scalable hyperbolic {PDE} engine}, Axioms, 7 (2018), p.~63,
  \url{https://doi.org/10.3390/axioms7030063}.

\bibitem{fehn2020hybrid}
{\sc N.~Fehn, P.~Munch, W.~A. Wall, and M.~Kronbichler}, {\em Hybrid multigrid
  methods for high-order discontinuous {G}alerkin discretizations}, J. Comput.
  Phys., 415 (2020), p.~109538,
  \url{https://doi.org/10.1016/j.jcp.2020.109538}.

\bibitem{Hackbusch85}
{\sc W.~Hackbusch}, {\em Multi-grid Methods and Applications}, Springer,
  Heidelberg, 1985.

\bibitem{hong2016robust}
{\sc Q.~Hong, J.~Kraus, J.~Xu, and L.~Zikatanov}, {\em A robust multigrid
  method for discontinuous {G}alerkin discretizations of {S}tokes and linear
  elasticity equations}, Numer. Math., 132 (2016), pp.~23--49.

\bibitem{JanssenKanschat11}
{\sc B.~Janssen and G.~Kanschat}, {\em Adaptive multilevel methods with local
  smoothing for {${H}^1$}- and {$H^{\text{curl}}$}-conforming high order finite
  element methods}, SIAM J. Sci. Comput., 33 (2011), pp.~2095--2114,
  \url{https://doi.org/10.1137/090778523}.

\bibitem{Kanschat08smoother}
{\sc G.~Kanschat}, {\em Robust smoothers for high order discontinuous
  {G}alerkin discretizations of advection-diffusion problems}, J. Comput. Appl.
  Math., 218 (2008), pp.~53--60,
  \url{https://doi.org/10.1016/j.cam.2007.04.032}.

\bibitem{KanschatMao15}
{\sc G.~Kanschat and Y.~Mao}, {\em Multigrid methods for {$\mathbf
  H^{\text{div}}$}-conforming discontinuous {G}alerkin methods for the {S}tokes
  equations}, J. Numer. Math., 23 (2015), pp.~51--66,
  \url{https://doi.org/10.1515/jnma-2015-0005}.

\bibitem{kohl2022textbook}
{\sc N.~Kohl and U.~R\"ude}, {\em Textbook efficiency: massively parallel
  matrix-free multigrid for the {S}tokes system}, SIAM J. Sci. Comput., 44
  (2022), pp.~C124--C155, \url{https://doi.org/10.1137/20M1376005}.

\bibitem{Kormann2011}
{\sc K.~Kormann and M.~Kronbichler}, {\em Parallel finite element operator
  application: {G}raph partitioning and coloring}, in Proceedings of the 7th
  IEEE International Conference on eScience, 2011, pp.~332--339,
  \url{https://doi.org/10.1109/eScience.2011.53}.

\bibitem{Kronbichler2012}
{\sc M.~Kronbichler and K.~Kormann}, {\em {A generic interface for parallel
  cell-based finite element operator application}}, Computers \& Fluids, 63
  (2012), pp.~135--147, \url{https://doi.org/10.1016/j.compﬂuid.2012.04.012}.

\bibitem{Kronbichler2017a}
{\sc M.~Kronbichler and K.~Kormann}, {\em Fast matrix-free evaluation of
  discontinuous {G}alerkin finite element operators}, ACM Trans. Math. Softw.,
  45 (2019), pp.~29/1--40, \url{https://doi.org/10.1145/3325864}.

\bibitem{kronbichler2019multigrid}
{\sc M.~Kronbichler and K.~Ljungkvist}, {\em Multigrid for matrix-free
  high-order finite element computations on graphics processors}, ACM Trans.
  Parallel Comput., 6 (2019), pp.~2/1--32,
  \url{https://doi.org/10.1145/3322813}.

\bibitem{kronbichler2022cg}
{\sc M.~Kronbichler, D.~Sashko, and P.~Munch}, {\em Enhancing data locality of
  the conjugate gradient method for high-order matrix-free finite-element
  implementations}, Int. J. High-Performance Computing Appl., 37 (2023),
  pp.~61--81, \url{https://doi.org/10.1177/10943420221107880}.

\bibitem{Lynch1964}
{\sc R.~E. Lynch, J.~R. Rice, and D.~H. Thomas}, {\em Direct solution of
  partial difference equations by tensor product methods}, Numer. Math., 6
  (1964), pp.~185--199, \url{https://doi.org/10.1007/BF01386067}.

\bibitem{malas2017multidimensional}
{\sc T.~M. Malas, G.~Hager, H.~Ltaief, and D.~E. Keyes}, {\em Multidimensional
  intratile parallelization for memory-starved stencil computations}, ACM
  Trans. Parallel Comput., 4 (2017), pp.~12:1--32,
  \url{https://doi.org/10.1145/3155290}.

\bibitem{Melenk99fullydiscrete}
{\sc J.~M. Melenk, K.~Gerdes, and C.~Schwab}, {\em Fully discrete hp-finite
  elements: fast quadrature}, Comp. Meth. Appl. Mech. Engrg., 190 (2001),
  pp.~4339--4364, \url{https://doi.org/10.1016/S0045-7825(00)00322-4}.

\bibitem{munch2023cache}
{\sc P.~Munch and M.~Kronbichler}, {\em Cache-optimized and low-overhead
  implementations of additive {S}chwarz methods for high-order {FEM} multigrid
  computations}, Int. J. High Perf. Comput. Appl.,  (2023),
  \url{https://doi.org/10.1177/10943420231217221}.
\newblock In press.

\bibitem{Orszag1980}
{\sc S.~A. Orszag}, {\em Spectral methods for problems in complex geometries},
  J. Comput. Phys., 37 (1980), pp.~70--92.

\bibitem{thomas2003textbook}
{\sc J.~L. Thomas, B.~Diskin, and A.~Brandt}, {\em Textbook multigrid
  efficiency for fluid simulations}, Annu. Rev. Fluid Mech., 35 (2003),
  pp.~317--340, \url{https://doi.org/10.1146/annurev.fluid.35.101101.161209}.

\bibitem{Trojak2022}
{\sc W.~Trojak, R.~Watson, and F.~Witherden}, {\em Hyperbolic diffusion in flux
  reconstruction: Optimisation through kernel fusion within tensor-product
  elements}, Computer Phys. Commun., 273 (2022), p.~108235,
  \url{https://doi.org/10.1016/j.cpc.2021.108235}.

\bibitem{WitteArndtKanschat21}
{\sc J.~Witte, D.~Arndt, and G.~Kanschat}, {\em Fast tensor product {S}chwarz
  smoothers for high-order discontinuous {G}alerkin methods}, Comput. Meth.
  Appl. Math., 21 (2021), pp.~709--728,
  \url{https://doi.org/10.1515/cmam-2020-0078}.

\end{thebibliography}
